\documentclass{article}

\usepackage[preprint]{neurips_2021}

\usepackage[utf8]{inputenc} 
\usepackage[T1]{fontenc}    
\usepackage{hyperref}       
\usepackage{url}            
\usepackage{booktabs}       
\usepackage{amsfonts}       
\usepackage{nicefrac}       
\usepackage{microtype}      

\usepackage[pdftex]{graphicx}
\usepackage{multirow}
\usepackage{adjustbox}
\usepackage{amsmath}
\usepackage{caption}

\title{Confidence Threshold Neural Diving}

\author{%
  Taehyun Yoon \\
  Artificial Intelligence Graduate School\\
  UNIST\\
  \texttt{thyoon@unist.ac.kr} \\
}

\begin{document}

\maketitle

\begin{abstract}
  Finding a better feasible solution in a shorter time is an integral part of solving Mixed Integer Programs. We present a post-hoc method based on Neural Diving [1] to build heuristics more flexibly. We hypothesize that variables with higher confidence scores are more definite to be included in the optimal solution. For our hypothesis, we provide empirical evidence that confidence threshold technique produces partial solutions leading to final solutions with better primal objective values. Our method won 2nd place in the primal task on the NeurIPS 2021 ML4CO competition. Also, our method shows the best score among other learning-based methods in the competition.  
\end{abstract}

\section{Introduction}
\label{intro}

  Mixed Integer Programming (MIP) solvers are dependent on various heuristics. It often requires parameter tuning or custom heuristics to solve problems with specific distributions for better optimization results. There have been approaches to partially replace existing rules and heuristics in MIP solvers. They tackle the components of MIP solvers with machine learning: node selection [3, 8], variable selection [4, 5, 6, 15], cutting-plane method [7, 9, 16] and primal heuristics [1, 10, 11, 12, 14]. In particular, we choose Neural Diving from [1] as a baseline to learn a diving-style primal heuristic.  
  
  The idea of Neural Diving is to learn to generate partial solutions given a problem represented in a bipartite graph. Formulated as a supervised learning problem, a graph neural network model is trained to generate a solution with a loss function to maximize the log likelihood with respect to the collected solutions. To generate solutions for partial discrete variables, [1] adopts SelectiveNet [2] to optimize for coverage among variables.  
  
  We find that the coverage control is critical for the performance in the competition. Meanwhile, there are two issues for optimizing the coverage with SelectiveNet. First, we need to set the coverage rate, which is unknown prior, when we train the model. It becomes problematic if the computation resource, especially GPU, is limited because we need to train the model for every coverage rate. Even if we train the model for every coverage rate, the actual coverage rate does not match precisely. Second, better validation loss and accuracy from using SelectiveNet does not lead to better performance in the competition. We do not know the actual validation performance of the model unless we measure the actual primal integral value with a time limit given in the competition.
  
  In this report, we suggest the confidence threshold technique as a better alternative to the SelectiveNet approach. The confidence threshold technique is free from the first issue since it is a deterministic post-hoc method. Also, we solve the second issue by performing as many evaluations to measure primal integral with different confidence levels as needed. We hypothesize that the confidence score of the model for each variable represents its level of definiteness to include the corresponding variable in the target solution. In Table \ref{tab:ML4CO}, we suggest empirical proof of our hypothesis by showing that the confidence threshold method is an effective solution in the competition.  

\begin{table*}[h!]
\centering
\caption{Leaderboard of the primal task on the NeurIPS 2021 ML4CO competition. Average over 100, 100 and 20 problem instances for item placement, load balancing and anonymous dataset are reported, respectively. Cumulative reward refers to the negated and unshifted version of the evaluation metric — primal integral.   \\} 
\label{tab:ML4CO}
  \begin{adjustbox}{max width=\textwidth}
\begin{tabular}{|l|l|c|c|c|}
\hline
                       &                                                      & item\_placement   & load\_balancing   & anonymous         \\ \cline{3-5} 
\multirow{-2}{*}{rank} & \multirow{-2}{*}{team}                               & cumulative reward & cumulative reward & cumulative reward \\ \hline
1                      & CUHKSZ\_ATD                                          & -3355.56          & -213467.31        & -45747908.15      \\ \hline
\textbf{2}             & {\textbf{UNIST-LIM-Lab (Ours)}} & -5902.27          & -214696.55        & -48690877.94      \\ \hline
3                      & MDO                                                  & -4798.50          & -216053.14        & -52610680.40      \\ \hline
\end{tabular}
\end{adjustbox}
\end{table*}

\section{Background}
\label{bg}

In the competition, we are given Mixed Integer Linear Programs (MILP) instances expressed as follows.  

\begin{equation}
\label{formulation}
\begin{aligned}
\underset{\mathbf{x}}{\arg \min ~} & \mathbf{c}^{\top} \mathbf{x} \\
\text { subject to } & \mathbf{A}^{\top} \mathbf{x} \leq \mathbf{b}, \\
& \mathbf{x} \in \mathbb{Z}^{p} \times \mathbb{R}^{n-p},
\end{aligned}
\end{equation}

where $\bold{c} \in \mathbb R^n$ represents linear objective coefficients. Variables are denoted by $\bold x$. $\bold{A} \in \mathbb R^{m \times n}$ and $\bold{b} \in \mathbb R^m$ denote coefficients and upper bounds of linear constraints, respectively. $n$ is the number of total variables, $m$ is the number of total constraints and $p \leq n$ is the number of integer variables.  

\begin{figure}[h!]
    \centering
	\includegraphics[scale=0.3, keepaspectratio]{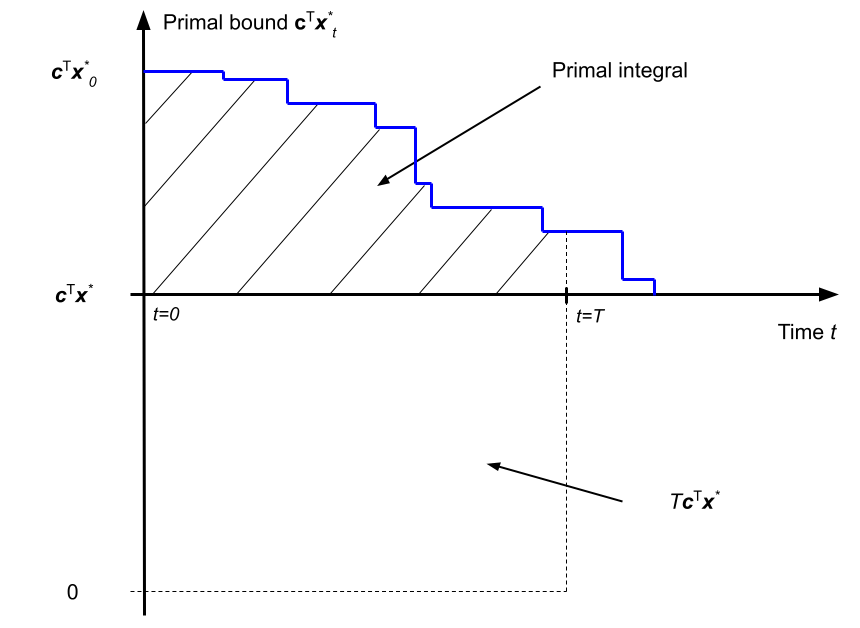}
	\caption{Illustration of primal bound and primal integral. This figure is an excerpt from the ML4CO competition description.}
	\label{fig:primal_integral}
\end{figure}  

\subsection{Primal task}
\label{primal}
In primal task, our goal is to minimize the average primal integral of test instances. Primal integral is defined as follows. 

\begin{equation}
\int_{t=0}^{T} \mathbf{c}^{\top} \mathbf{x}_{t}^{\star} \mathrm{d} t-T \mathbf{c}^{\top} \mathbf{x}^{\star}
\end{equation}

where $\bold{x}_t^\star$ is the best feasible solution found so far at time $t$, which makes $\mathbf{c}^{\top} \mathbf{x}_{t}^{\star}$ primal bound at time $t$. $T$ is a time limit for each instance. $T \mathbf{c}^{\top} \mathbf{x}^{\star}$ is an instance-specific constant that depends on the optimal solution $\mathbf{x}^{\star}$. To illustrate, primal integral is the area in slashes in Figure \ref{fig:primal_integral}.

\section{Method}
\label{method}

\begin{figure}[h!]
    \centering
	\includegraphics[width=\textwidth, keepaspectratio]{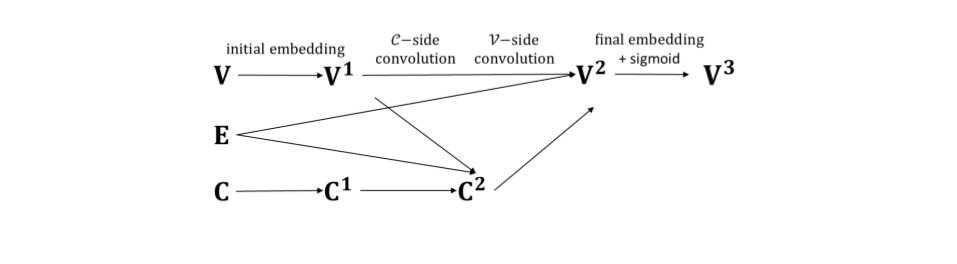}
	\caption{GCNN from learn2branch [4].}
	\label{fig:gnn}
\end{figure}  

\begin{figure}[htp!]
    \centering
    \captionsetup{width=0.45\linewidth}
    \begin{minipage}[t]{0.5\linewidth}
	\includegraphics[width=\linewidth, keepaspectratio]{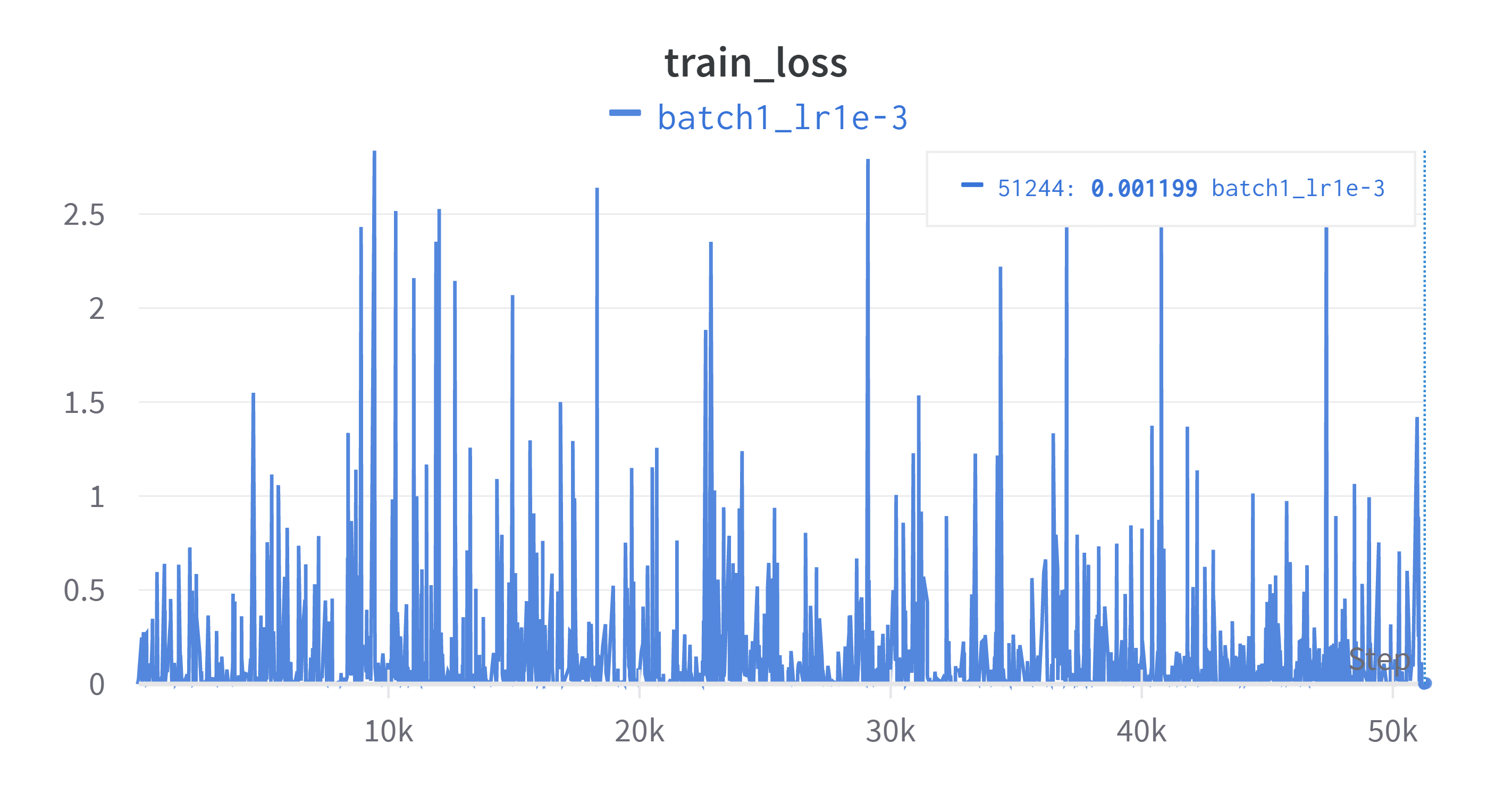}
	\caption{Training loss plot for the anonymous dataset with the full-batch loss function.}
	\label{fig:train_fullbatch_loss}
    \end{minipage}%
    \hfill%
    \begin{minipage}[t]{0.5\linewidth}
	\includegraphics[width=\linewidth, keepaspectratio]{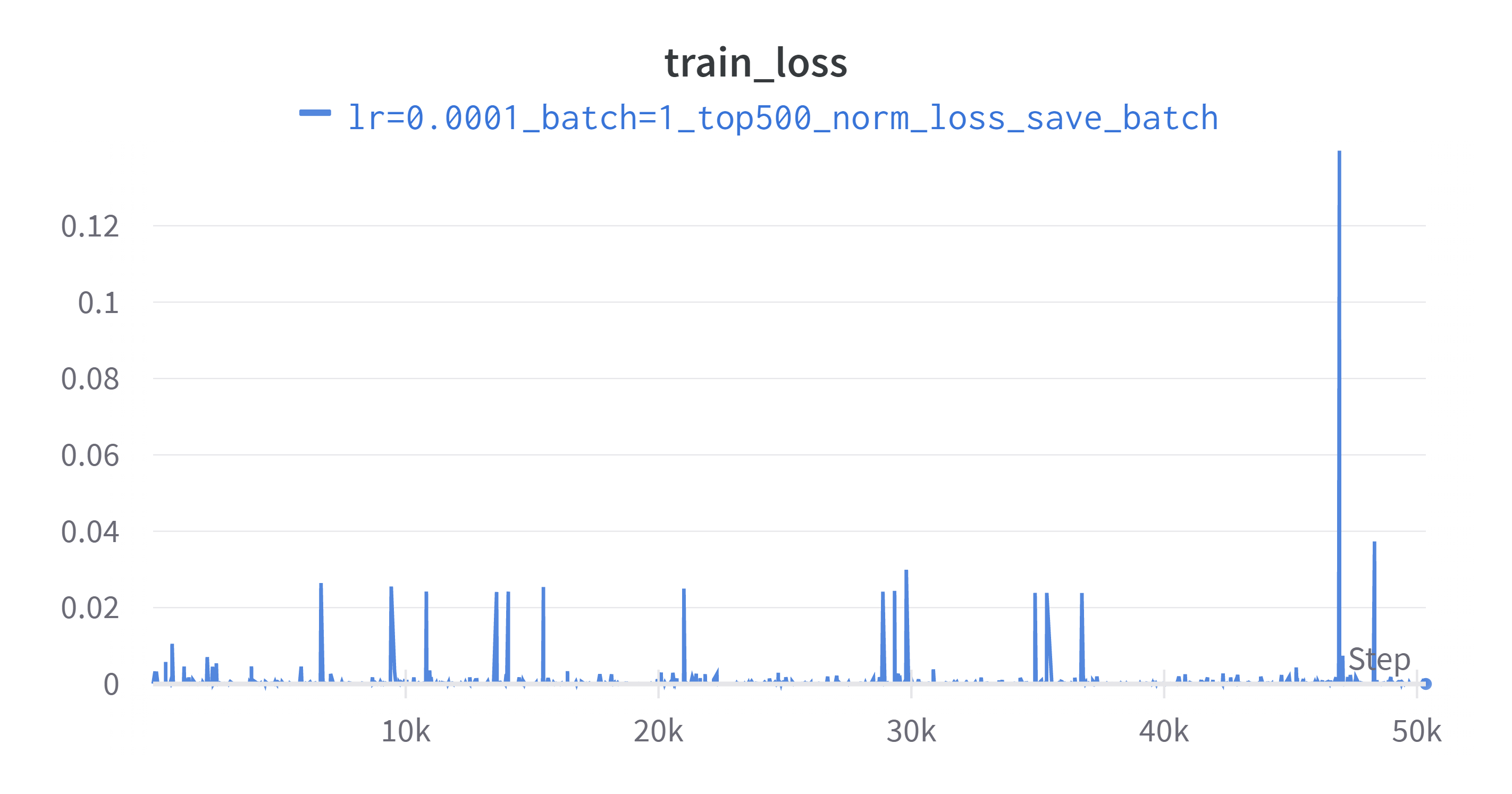}
	\caption{Training loss plot for the anonymous dataset with the mini-batch loss function.}
	\label{fig:train_minibatch_loss}
    \end{minipage}%
\end{figure}

Our method is based on Neural Diving from [1]. On the other hand, the differences in our method from [1] are as follows. First, we use GCNN formulated in [4] with two half convolutions as shown in Figure \ref{fig:gnn}. Second, we suggest the confidence threshold component in place of the SelectiveNet [2] component. Third, we construct a mini-batch loss function that normalizes the loss by the number of total nodes. $N_B$ is the number of graphs in a mini-batch and $N_i$ is the number of nodes in a graph. The reason why we construct the mini-batch loss is that the full batch loss from [1] shows unstable training loss for the anonymous dataset as the instances show a big difference in distribution. Figure \ref{fig:train_fullbatch_loss} and \ref{fig:train_minibatch_loss} shows that mini-batch loss stabilizes the training the GCNN model for the anonymous dataset.  

\begin{equation}
L_{\text {minibatch }}=-\frac{1}{N_{B}} \sum_{i=1}^{N_{B}} \frac{1}{N_{i}} \sum_{j=1}^{N_{i}} w_{i, j} \log p_{\theta}\left(x^{i, j} \mid M_{i}\right)
\end{equation}

\begin{figure}[htp!]
    \centering
	\includegraphics[width=\textwidth, keepaspectratio]{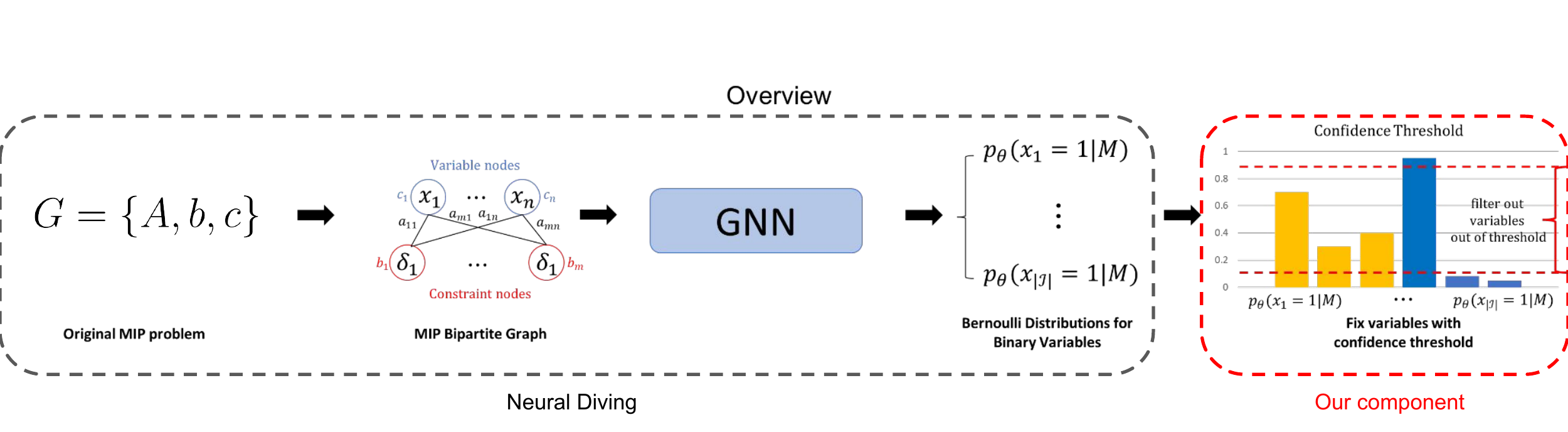}
	\caption{The overview of Confidence Threshold Neural Diving.}
	\label{fig:overview}
\end{figure}

The overview of our method is illustrated in Figure \ref{fig:overview}.   

\subsection{Solution collection}
We train neural network models for each dataset given from the competition separately. First, we collect solutions of the MIP problems with the solver in a transformed problem space. For each dataset, item placement, load balancing and anonymous, we give time limits of 2 hours, 20 hours and 20 hours, respectively. We set heuristic emphasis to aggressive  in common without tuning specific parameters.  

\subsubsection{Data augmentation}
For anonymous dataset, we are given only 100 training instances. To handle this data scarcity issue, we augment the dataset with MIPLIB 2010 collection problems, where we obtain 599 {graph sample, solution} pairs additionally.  

\subsection{SCIP parameter tuning}
For item placement dataset, we do not use Neural Diving model, but only tune the SCIP parameters. We suppose the neural network model was unable to capture the underlying pattern of the given multi-dimensional knapsack problems whose objective coefficients of discrete variables are all 0 and the ratio of solution variables assigned 1 is fixed. Instead, we discover the following SCIP parameters to improve the primal integral score on the item placement dataset: ‘lp/solvedepth=24’ and ‘separating/maxcuts=0’. For other datasets, load balancing and anonymous, we do not tune the parameters of the solver at test time but we set heuristic emphasis to aggressive.  

\subsection{Evaluation}
We unfix the variables and switch to use the solver only if the model-generated partial solution turns out to be infeasible. We evaluate each model for load balancing and anonymous dataset by grid searching of the symmetric confidence thresholds as shown in the Figure \ref{fig:overview}.

\section{Results}
\label{results}

\begin{figure}[htp!]
    \centering
    \captionsetup{width=0.45\linewidth}
    \begin{minipage}[t]{0.5\linewidth}
	\includegraphics[width=\linewidth, keepaspectratio]{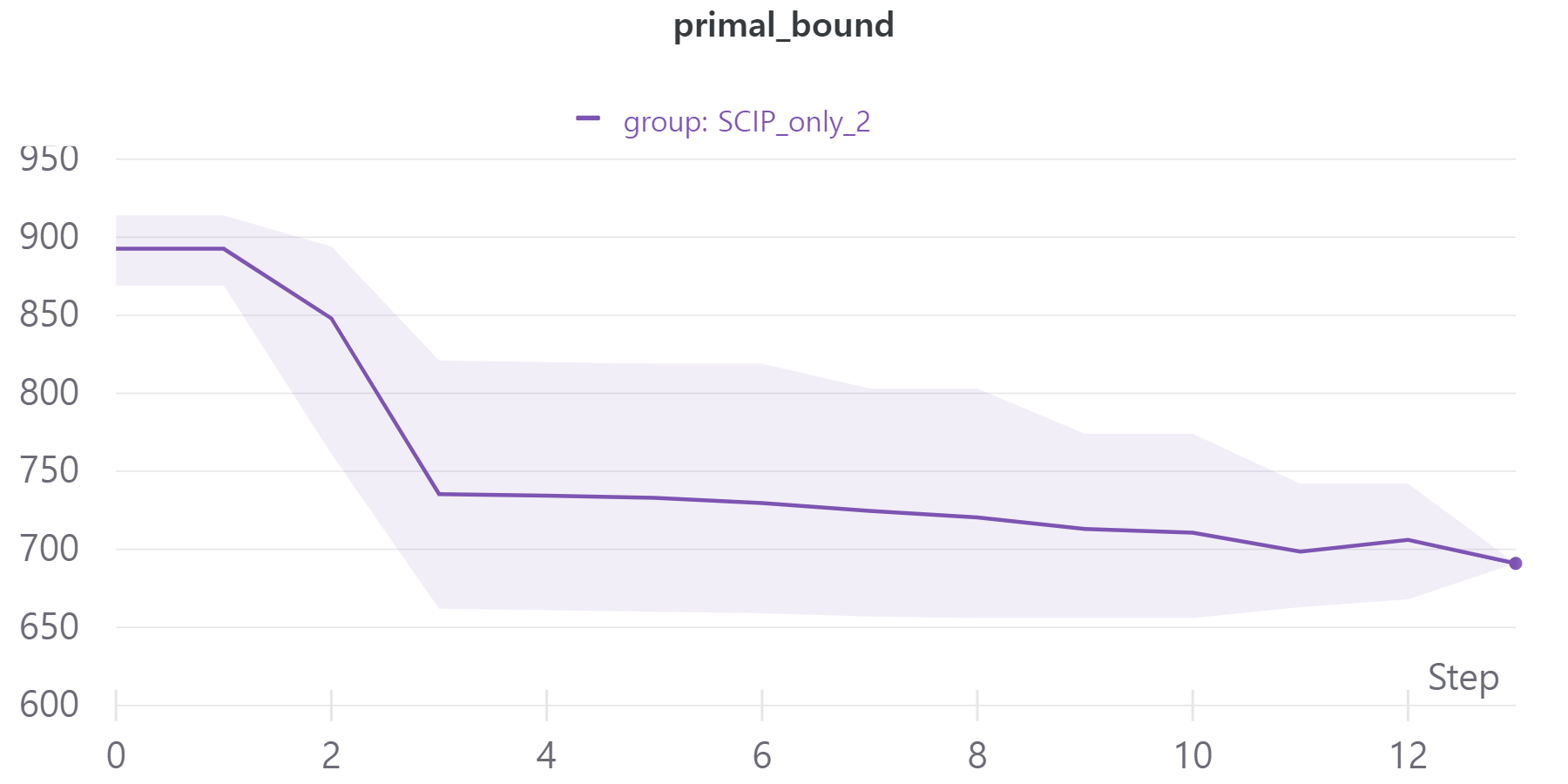}
	\caption{Primal bound of the SCIP solver with aggressive heuristic emphasis over optimization steps for the load balancing dataset.  }
	\label{fig:SCIP_primal_bound}
    \end{minipage}%
    \hfill%
    \begin{minipage}[t]{0.5\linewidth}
	\includegraphics[width=\linewidth, keepaspectratio]{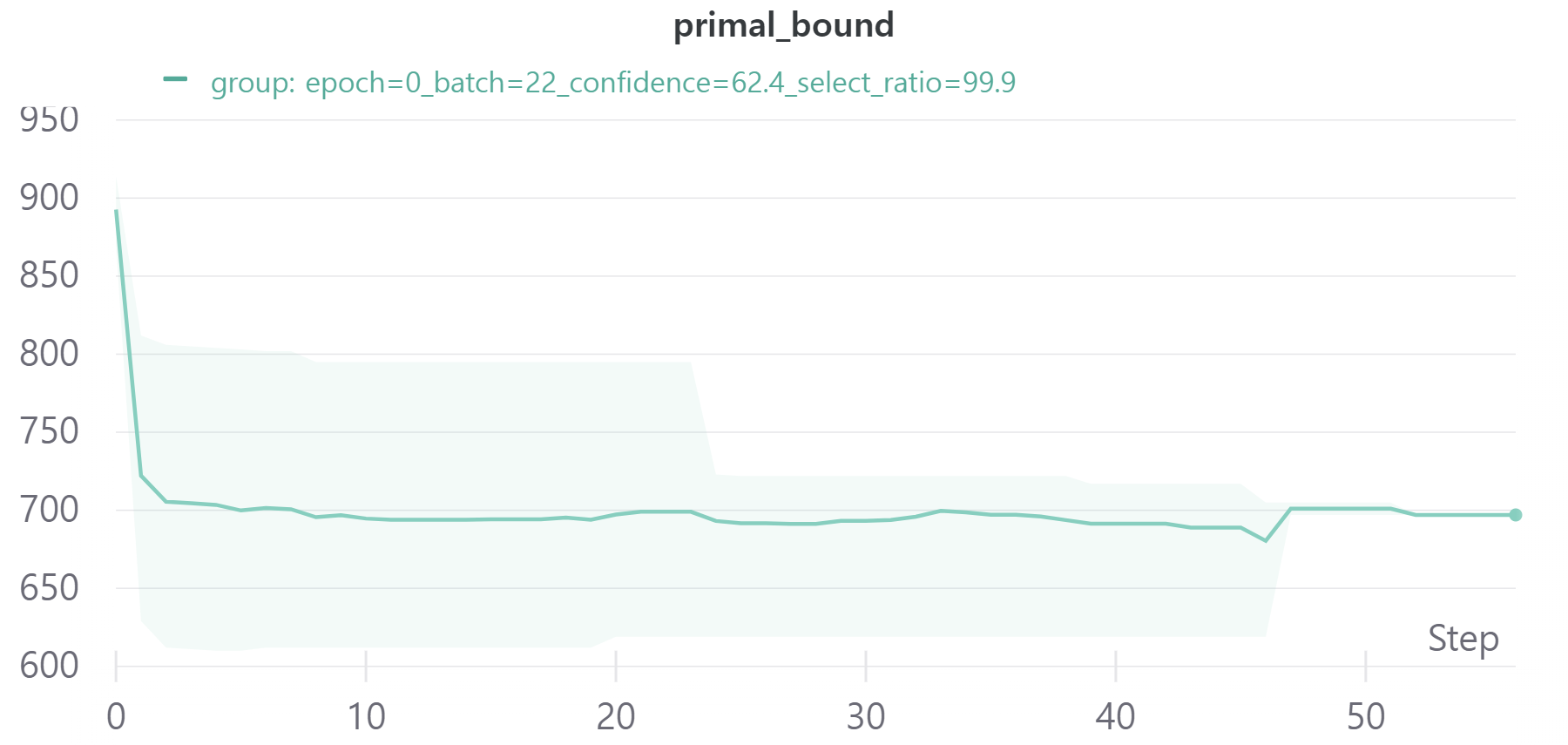}
	\caption{Primal bound of Confidence Threshold Neural Diving over optimization steps for the load balancing dataset.  }
	\label{fig:CTND_primal_bound}
    \end{minipage}%
\end{figure}

We find that it shows the best validation primal integral score when it fixes 96.8\% variables by confidence threshold of 62.4\% for the load balancing dataset, showing 100\% feasibility. On the other hand, for anonymous dataset, we observe the best validation score when it fixes 70.39\% of variables by confidence threshold of 93\%, showing 40\% feasibility. Figure \ref{fig:SCIP_primal_bound} and \ref{fig:CTND_primal_bound} show the primal bounds of the SCIP solver with aggressive heuristic emphasis and Confidence Threshold Neural Diving, respectively, over optimization steps. Apparently, Confidence Threshold Neural Diving drops the primal bound much faster than the SCIP solver with aggressive heuristic emphasis. Finally, the test results in the competition in Table \ref{tab:ML4CO} serve as empirical evidence that Confidence Threshold Neural Diving is a competitive method.    

\subsection{Conclusion}

We show that Neural Diving combined with the confidence threshold technique is a competitive method to find better primal solutions in a shorter time. The competition result empirically supports our hypothesis that variables with higher confidence scores are more definite to be included in the target solution with a better objective value. By using post-hoc method, confidence threshold, we avoid the issue of training the model for every coverage rate. Also, we resolve the discrepancy between the ML objective (loss and accuracy) and the MIP objective (primal integral) by performing greed search to find the coverage rate and the weights that result in better primal integral value.

\begin{ack}
This work was supported by the 2021 Research Fund (1.210107.01) of UNIST and National Research Foundation of Korea(NRF) funded by the Korea government(MSIT)(2021R1A4A3033149).
We would like to thank Minsub Lee, Keehun Park, Hanbum Ko, Jinwon Choi, Hyunho Lee, Hyokun Yun and Sungbin Lim for their advice and constructive discussions. Finally, we appreciate Kakao Enterprise for sharing Brain Cloud during the competition.  

\end{ack}

\section*{References}

{
\small

[1] Nair, V., \ Bartunov, S., \ Gimeno, F., \ von Glehn, I., \ Lichocki, P., \ Lobov, I., ... \& Zwols, Y. (2020). Solving mixed integer programs using neural networks. {\it arXiv preprint arXiv:2012.13349.}   

[2] Geifman, Y., \& El-Yaniv, R. (2019, May). Selectivenet: A deep neural network with an integrated reject option. {\it In International Conference on Machine Learning (pp. 2151-2159). PMLR.}  

[3] He, H., Daume III, H., \& Eisner, J. M. (2014). Learning to search in branch and bound algorithms. {\it Advances in neural information processing systems, 27, 3293-3301.}  

[4] Gasse, M., Chételat, D., Ferroni, N., Charlin, L., \& Lodi, A. (2019). Exact combinatorial optimization with graph convolutional neural networks. {\it arXiv preprint arXiv:1906.01629.}  

[5] Gupta, P., Gasse, M., Khalil, E. B., Kumar, M. P., Lodi, A., \& Bengio, Y. (2020). Hybrid models for learning to branch. {\it arXiv preprint arXiv:2006.15212.}  

[6] Sun, H., Chen, W., Li, H., \& Song, L. (2020). Improving learning to branch via reinforcement learning.  

[7] Tang, Y., Agrawal, S., \& Faenza, Y. (2020, November). Reinforcement learning for integer programming: Learning to cut. {\it In International Conference on Machine Learning (pp. 9367-9376). PMLR.}  

[8] Yilmaz, K., \& Yorke-Smith, N. (2020). Learning efficient search approximation in mixed integer branch and bound. {\it arXiv preprint arXiv:2007.03948.}  

[9] Ding, J. Y., Zhang, C., Shen, L., Li, S., Wang, B., Xu, Y., \& Song, L. (2020, April). Accelerating primal solution findings for mixed integer programs based on solution prediction. {\it  In Proceedings of the AAAI Conference on Artificial Intelligence (Vol. 34, No. 02, pp. 1452-1459).}  

[10] Song, J., Lanka, R., Yue, Y., \& Dilkina, B. (2020). A general large neighborhood search framework for solving integer linear programs. {\it arXiv preprint arXiv:2004.00422.}  

[11] Shen, Y., Sun, Y., Eberhard, A., \& Li, X. (2021, July). Learning Primal Heuristics for Mixed Integer Programs. {\it In 2021 International Joint Conference on Neural Networks (IJCNN) (pp. 1-8). IEEE.}  

[12] Sonnerat, N., Wang, P., Ktena, I., Bartunov, S., \& Nair, V. (2021). Learning a large neighborhood search algorithm for mixed integer programs. {\it arXiv preprint arXiv:2107.10201.}  

[13] Paulus, A., Rolínek, M., Musil, V., Amos, B., \& Martius, G. (2021). CombOptNet: Fit the Right NP-Hard Problem by Learning Integer Programming Constraints. {\it arXiv preprint arXiv:2105.02343.}  

[14] Qi, M., Wang, M., \& Shen, Z. J. (2021). Smart Feasibility Pump: Reinforcement Learning for (Mixed) Integer Programming. {\it arXiv preprint arXiv:2102.09663.}

[15] Zarpellon, G., Jo, J., Lodi, A., \& Bengio, Y. (2020). Parameterizing branch-and-bound search trees to learn branching policies. {\it arXiv preprint arXiv:2002.05120, 12.}

[16] Huang, Z., Wang, K., Liu, F., Zhen, H. L., Zhang, W., Yuan, M., ... \& Wang, J. (2022). Learning to select cuts for efficient mixed-integer programming. {\it Pattern Recognition, 123, 108353.}

}

\end{document}